\documentclass[twocolumn]{autart}    

\usepackage{graphicx}          
\usepackage{color}                                

\begin{document}

\begin{frontmatter}

\title{Regional controllability analysis of fractional diffusion equations
with Riemann-Liouville time fractional derivatives} 

\thanks[footnoteinfo]{
This work was completed while F. Ge visited the Mechatronics, Embedded Systems and Automation Lab,
 University of California, Merced from October 2014 to October
2015.\\\\
  $^{\dag}$Corresponding author. Tel./Fax:+1(209)228-4672/4047.}

\author[CUG1,CUG2]{Fudong Ge}\ead{gefd2011@gmail.com},    
\author[UCMerced]{YangQuan Chen}$^{,\dag}$
\ead{yqchen@ieee.org},   
\author[DHU2]{Chunhai Kou}\ead{kouchunhai@dhu.edu.cn}  

\address[CUG1]{School of Computer Science, China University of Geosciences, Wuhan 430074, PR China}
\address[CUG2]{Hubei Key Laboratory of Intelligent Geo-Information Processing, China University of Geosciences, Wuhan 430074, PR China} 
\address[UCMerced]{Mechatronics, Embedded Systems and Automation Lab,
 University of California, Merced, CA 95343, USA }             
\address[DHU2]{Department of Applied Mathematics,
 Donghua University, Shanghai 201620, PR China}        


\begin{abstract}                          
This paper is concerned with the concepts of regional controllability for the Riemann-Liouville
time fractional diffusion systems of order $\alpha\in(0,1)$. The characterizations of strategic
actuators to achieve regional controllability are investigated
when the control inputs emerge in the differential equations as distributed inputs.  In the end,
an approach to guarantee the regional controllability of the problems
under consideration in the considered subregion with minimum energy control is described and successfully
tested through two applications.
\end{abstract}

\end{frontmatter}

\section{Introduction}

Recently sub-diffusion processes have attracted increasing interest since the introduction of continuous time random walks (CTRWs) in \cite{Montroll1965} and a large number of contributions have been given to them (\cite{mainardi2007sub,ralf},
\cite{Ge2016JMAA,japankf1}). Since CTRW is a random walk subordinated to a simple renewal process, by \cite{hilfer1995fractional},
it can be regarded as a
generalized physical diffusion process (including the sub-diffusion process and the super-diffusion process)
and there exists a closed connection between  the time fractional diffusion system and the sub-diffusion process.
Moreover, it is confirmed in \cite{ralf} and  \cite{Mandelbrot} that the time fractional diffusion systems can be used to  well
characterize those sub-diffusion processes, which offer better performance not achievable before using conventional diffusion systems and
surely raise many potential research opportunities at the same time.


In the case of diffusion system, it is  well known that in general, not all the states can be reached
in the whole domain of interest. So here, we first introduce some notations
on the regional controllability of  time fractional diffusion systems  when the system under consideration is
only exactly (or approximately) controllability  on a subset of the whole
space,  which can be regarded as an extensions  of the research work in
(\cite{214}, \cite{179}). Besides, focusing on regional controllability would allow
for a reduction in the number of physical actuators,  offer the potential to reduce computational requirements in some cases, and  also possible to
discuss those systems which are not controllable on the whole domain, etc.

Furthermore, in  \cite{chen2016adaptive,Ge2016IJC} and \cite{AEIJAI}, the authors have shown that the measurements and actions in practical systems
can be better described by using the notion of actuators  and sensors (including the location,
number and spatial distribution of actuators and sensors \cite{AEIJAI2}).
Then the contribution of this present work is on the regional controllability of the sub-diffusion processes described by Riemann-Liouville time fractional diffusion systems of order $\alpha\in (0,1)$ by using the notion of actuators  and sensors.
As cited
in \cite{hilfer2000applications},  their applications
are rich in many real life. For example, the flow through porous media (\cite{uchaikin2012fractional}),
or the swarm of robots moving  through dense forest (\cite{2Spears}). We hope that the results
here could provide some insights into the qualitative analysis of the
design and configuration of fractional controller.

The rest of the paper is organized as follows. The mathematical concept of regional
controllability problem is presented in the next section. Section 3 is focused on the
characterizations of strategic actuators in the case of regional controllability.
In Section 4, our main results on the regional controllability analysis of time time fractional diffusion systems
are presented and the determination of the optimal control which achieves the regional
controllability is obtained. Two applications are worked out in the last section.
\section{Statement of the problem}\label{sec:2}

Let $\Omega$ be an open bounded subset of $\mathbf{R}^n$ with smooth boundary $\partial\Omega$ and we
consider the following abstract Riemann-Liouville time fractional differential system:
\begin{equation}\label{problem}
\left\{\begin{array}{l}
_0D^{\alpha}_tz(t)=Az(t)+Bu(t), ~~ t\in [0,b],~0<\alpha<1,
\\ \lim\limits_{t\to 0^+} {}_0D^{\alpha}_t z(t)=z_0,
\end{array}\right.
\end{equation}
where $A$ generates a strongly continuous semigroup $\{\Phi(t)\}_{t\geq 0}$ on the Hilbert
space $Z:=L^2(\Omega)$, $-A$ is a uniformly elliptic operator (\cite{uniformlyelliptic}, \cite{uniformlyelliptic2}),
$z\in L^2(0,b;Z)$ and the initial vector $z_0\in Z$.
Here
$_0D^{\alpha}_t$ and $_0I^{\alpha}_t$ denote the Riemann-Liouville
fractional order derivative and integral, respectively, given by \cite{Kilbas}
\begin{equation}\label{caputoderivative}
_0D_t^\alpha z(t)=\frac{d}{d t}{}_0I_t^{1-\alpha}z(t),~~0<\alpha<1
\end{equation}
and
\begin{equation}{}_0I^{\alpha}_tz(t)=\frac{1}{\Gamma(\alpha)}\int^{t}_{0}{(t-s)^{\alpha-1}z(s)ds},~~\alpha>0.\end{equation}
In addition,  $B$ is a control operator depends on the number and the structure of actuators. The control $u\in U$
where $U$ is a Hilbert space. In particular, if the system is excited by $p$ actuators, one has $u\in L^2(0,b;\mathbf{R}^p)$
and $B: \mathbf{R}^p \to Z.$

We first recall some necessary lemmas to be used afterwards.
 \begin{lem}\label{Lem0}
 For any given $f\in  L^2\left(0,b; Z\right),$ $0<\alpha<1,$
a function $v\in L^2\left(0,b; Z\right)$  is said to be a mild solution of the following system
\begin{equation}\label{problem2}
\left\{\begin{array}{l}
{}_0D^{\alpha}_tv(t)=Av(t)+f(t),~t\in [0,b],
\\ \lim\limits_{t\to 0^+}{}_0D_t^{\alpha-1}v(t)=v_0\in Z,
\end{array}\right.
\end{equation}
if it satisfies
\begin{eqnarray}\label{transform}
v(t)=t^{\alpha-1}K_\alpha(t)v_0+\int_0^t(t-s)^{\alpha-1}{K_\alpha(t-s)f(s)}ds,
\end{eqnarray}
where
\begin{eqnarray}\label{Salpha}
K_\alpha(t)=\alpha\int_0^\infty{\theta\phi_\alpha(\theta)\Phi(t^\alpha\theta)}d\theta.\end{eqnarray}
Here $\{\Phi(t)\}_{t\geq 0}$ is the strongly continuous semigroup generated by operator $A$, $\phi_\alpha(\theta)=
\frac{1}{\alpha}\theta^{-1-\frac{1}{\alpha}}\psi_\alpha(\theta^{-\frac{1}{\alpha}})$ and  $\psi_\alpha$ is a probability density function defined by
$( \theta >0)$
\begin{eqnarray}\label{pf}
\psi_\alpha(\theta)=\frac{1}{\pi}\sum\limits_{n = 1}^\infty {( - 1)^{n - 1}\theta^{-\alpha n-1}\frac{\Gamma(n\alpha+1)}{n!}
\sin(n\pi \alpha)}\end{eqnarray}
such that (\cite{Mainardi})
\begin{eqnarray*}\int_0^\infty{\psi_\alpha(\theta)}d\theta=1\mbox{ and }\int_0^\infty{\theta^\nu
\phi_\alpha(\theta)}d\theta=\frac{\Gamma(1+\nu)}{\Gamma(1+\alpha\nu)},~\nu\geq 0.\end{eqnarray*}
 \end{lem}
\textbf{Proof.}
It follows from the Laplace transforms
 \begin{eqnarray}\tilde{v}(\lambda)=\int_0^\infty{e^{-\lambda s}v(s)}ds~\mbox{ and }~\tilde{f}(\lambda)=\int_0^\infty{e^{-\lambda s}f(s)}ds\end{eqnarray}
 that the system $(\ref{problem})$ is equivalent to (\cite{lin2013laplace})
 \begin{eqnarray}
~\lambda^\alpha \tilde{v}(\lambda)-v_0-A\tilde{v}(\lambda)=\tilde{f}(\lambda).
\end{eqnarray}
Then
 \begin{equation}
\begin{array}{l}
\tilde{v}(\lambda)=(\lambda^\alpha I-A)^{-1}(v_0+\tilde{f}(\lambda))\\{\kern 18pt}
=\int_0^\infty{e^{-\lambda^\alpha s}\Phi(s)[v_0+\tilde{f}(\lambda)]}ds.
\end{array}
\end{equation}
Consider the stable probability density function\index{probability density function} $(\ref{pf})$. By the arguments in \cite{Mainardi}, we see that $\psi_\alpha(\theta)$$(\theta>0)$
satisfies the following property
 \begin{eqnarray}\widetilde{\psi_\alpha}(\lambda)=\int_0^\infty{e^{-\lambda \theta}\psi_\alpha(\theta)}d\theta=e^{-\lambda^\alpha},~~ \alpha\in (0,1).\end{eqnarray}
Let $s=\tau^\alpha.$  We obtain that
 \begin{eqnarray*}
 \tilde{v}(\lambda)&=&
 \alpha\int_0^\infty{e^{-\lambda^\alpha \tau^\alpha}\Phi(\tau^\alpha)\tau^{\alpha-1}[v_0+\tilde{f}(\lambda)]}d\tau\\&=&
 \alpha\int_0^\infty{\int_0^\infty{e^{-\lambda \tau\theta}\psi_\alpha(\theta)\Phi(\tau^\alpha)\tau^{\alpha-1}[v_0+\tilde{f}(\lambda)]}d\theta}d\tau\\&=&
 \sigma_1(v_0)+\sigma_2(f),
\end{eqnarray*}
where $
 \sigma_1(v_0)= \alpha\int_0^\infty{\int_0^\infty{e^{-\lambda \tau\theta}\psi_\alpha(\theta)\Phi(\tau^\alpha)\tau^{\alpha-1}}d\theta}d\tau v_0
$
 and
 $
\sigma_2(f)= \alpha\int_0^\infty{\int_0^\infty{e^{-\lambda \tau\theta}\psi_\alpha(\theta)\Phi(\tau^\alpha)\tau^{\alpha-1}\tilde{f}(\lambda)}d\theta}d\tau.
$
Suppose that $t=\tau \theta$. Then we have
 \begin{eqnarray*}
\begin{array}{l}\sigma_1(v_0)=
 \alpha\int_0^\infty{\int_0^\infty{e^{-\lambda t}\psi_\alpha(\theta)\Phi\left(\frac{t^\alpha}
 {\theta^\alpha}\right)\frac{t^{\alpha-1}}{\theta^\alpha}}d\theta}dt
  ~v_0\\=\int_0^\infty{e^{-\lambda t}\alpha\int_0^\infty{\psi_\alpha(\theta)\Phi
  \left(t^\alpha\theta^{-\alpha}\right)t^{\alpha-1}\theta^{-\alpha}}d\theta}dt ~v_0\\
  =\int_0^\infty{e^{-\lambda t} t^{\alpha-1}\alpha\int_0^\infty{\frac{1}{\alpha}\theta^{-1-\frac{1}{\alpha}}\psi_\alpha(\theta^{-\frac{1}{\alpha}})
\theta\Phi\left(t^\alpha\theta\right)}d\theta}dt ~v_0
\end{array}
\end{eqnarray*}
and
 \begin{eqnarray*}
\begin{array}{l}
\sigma_2(f)\\
=\alpha\int_0^\infty{\int_0^\infty{\int_0^\infty{e^{-\lambda \tau\theta}\psi_\alpha(\theta)\Phi(\tau^\alpha)\tau^{\alpha-1}
e^{-\lambda s}f(s)}ds}d\theta}d\tau\\=
\alpha\int_0^\infty{\int_0^\infty{\int_0^\infty{e^{-\lambda (t+s)}\psi_\alpha(\theta)\Phi
\left(\frac{t^\alpha}{\theta^\alpha}\right)\frac{t^{\alpha-1}}{\theta^\alpha}f(s)}ds}d\theta}dt
\\=
\int_0^\infty{e^{-\lambda t}\alpha\int_0^t{\int_0^\infty{\psi_\alpha(\theta)\Phi
\left(\frac{(t-s)^\alpha}{\theta^{\alpha}}\right)\frac{(t-s)^{\alpha-1}f(s)}{\theta^{\alpha}}}d\theta}ds}dt
\\=
\int_0^\infty{e^{-\lambda t}\alpha\int_0^t{\int_0^\infty{\theta\frac{1}{\alpha}\theta^{-1-\frac{1}{\alpha}}\psi_\alpha(\theta^{-\frac{1}{\alpha}})
\frac{\Phi\left((t-s)^\alpha\theta\right)f(s)}{(t-s)^{1-\alpha}}}d\theta}ds}dt.
\end{array}
\end{eqnarray*}
Let $\phi_\alpha(\theta)=
\frac{1}{\alpha}\theta^{-1-\frac{1}{\alpha}}\psi_\alpha(\theta^{-\frac{1}{\alpha}})$ and \\
$K_\alpha(t)=
\alpha\int_0^\infty{\theta\phi_\alpha(\theta)\Phi(t^\alpha\theta)}d\theta.$
Then we get
\begin{eqnarray*}
v(t)=t^{\alpha-1}K_\alpha(t)v_0+\int_0^t{(t-s)^{\alpha-1}K_\alpha(t-s)f(s)}ds
\end{eqnarray*}
and the proof is complete.

\begin{lem}\cite{Bernard}\label{Lem2}
Let $\Omega \subseteq \mathbf{R}^n$ be
an open set and $C_0^{\infty}(\Omega)$ be the class of  infinitely differentiable functions on $\Omega$
 with compact support in $\Omega$ and $u\in L^1_{loc}(\Omega)$ be such that
\begin{eqnarray}
~~~~~~~~~~~~~~\int_\Omega u(x)\psi(x)dx=0,~~~\forall\psi \in C_0^{\infty}(\Omega).
\end{eqnarray}
Then $u=0$ almost everywhere in $\Omega.$
\end{lem}
\begin{lem}\cite{Klimek}\label{Lem1}
 Let the reflection operator $Q$ on interval $[0, b]$ be as follows:
\begin{eqnarray}\label{Q}
Qf(t) := f( b - t).\end{eqnarray}
Then the following equations hold:
\begin{eqnarray}\label{2.6}
Q_0I_t^\alpha f(t)={}_tI_b^\alpha Qf(t),~~~Q{}_0D_t^\alpha f(t)={}_tD_b^\alpha Qf(t)\end{eqnarray}
and
\begin{eqnarray}\label{2.7}
_0I_t^\alpha Qf(t) = Q{}_tI_b^\alpha f(t),~~~{}_0D_t^\alpha Qf(t)=Q{}_tD_b^\alpha f(t) .\end{eqnarray}
\end{lem}

Let $\omega \subseteq \Omega $ be a given region of positive Lebesgue measure and
$z_b\in L^2(\omega)$$($the target function$)$ be a given element. By Lemma~\ref{Lem0}, the unique mild solution $z(.,u)$ of $(\ref{problem})$ can be  given by
\begin{eqnarray*}
z(t,u)=t^{\alpha-1}K_\alpha(t)z_0+\int_0^t(t-s)^{\alpha-1}{K_\alpha(t-s)Bu(s)}ds.
\end{eqnarray*}
Taking into  account that $(\ref{problem})$ is a line system,
by the Proposition 3.1 in \cite{Ge2016IJC}, it suffices to suppose that $z_0=0$ in the following discussion.  Let $H:L^2(0,b;\mathbf{R}^p) \to  Z$  be
 \begin{eqnarray}\label{H}
Hu=\int_0^b{\frac{K_\alpha(b-s)}{(b-s)^{1-\alpha}}Bu(s)}ds,~\forall u\in L^2(0,b;\mathbf{R}^p)_.
\end{eqnarray}
In order to state the main results, the following two assumptions are supposed to hold all over the article:

$(A_1)$  $B$ is a densely defined operator and  $B^*$ exists.

$(A_2)$ $(BK_\alpha(t))^*$ exists and $(BK_\alpha(t))^*=K_\alpha^*(t)B^*$.

In particular, when $B\in \mathcal{L}\left(\mathbf{R}^p, Z\right)$ is a bounded linear operator from
$\mathbf{R}^p$ to $Z$, it is easy to see that $(A_1)$ and $(A_2)$ hold.
Suppose that $\{\Phi^*(t)\}_{t\geq 0}$, generated by the adjoint operator of $A$, is also a strongly continuous
semigroup in the space $Z$.
For any $v\in L^2(\Omega),$ by $\left<Hu,v\right>=\left<u,H^*v\right>$, we have
\begin{eqnarray}\label{Hstar}
H^* v=B^*(b-s)^{\alpha-1}K_\alpha^*(b-s)v,
\end{eqnarray} where $\left<\cdot,\cdot\right>$ is
the duality pairing of the space $Z$, $B^*$ is the adjoint operator of $B$ and
$K^*_\alpha(t)= \alpha\int_0^\infty{\theta\phi_\alpha(\theta)\Phi^*(t^\alpha\theta)}d\theta.$
Consider now the restriction map
\begin{eqnarray}p_\omega: L^2(\Omega)\to L^2(\omega),\end{eqnarray}
defined by $p_\omega z=z|_\omega$, is the projection operator on $\omega$. Then the adjoint operator of $p_\omega$ can be given by
\begin{eqnarray}\label{pomega}
p_\omega^* z(x):=\left\{\begin{array}{l}z(x), ~~x\in \omega ,  \\0,~~~x\in \Omega\backslash \omega.\end{array}\right.
\end{eqnarray}
and we are ready to state the following definition.

\begin{defn} \label{Def2}
$(i)$ The system $(\ref{problem})$ is said to be regionally exactly controllable on $\omega$
 if for any $z_b\in L^2(\omega)$ \textcolor[rgb]{1.00,0.00,0.00}{at time $b$}, there exists a control $u\in L^2(0,b; \mathbf{R}^p)$  such that
\begin{eqnarray}
p_\omega z(b,u)=z_b.
\end{eqnarray}
 $(ii)$ The system $(\ref{problem})$ is said to be regionally approximately controllable on
 $\omega$ at time $b$ if for any $z_b\in L^2(\omega)$, given $\varepsilon>0,$
 there exists a control $u\in L^2(0,b; \mathbf{R}^p)$  such that
\begin{eqnarray}
\|p_\omega z(b,u)-z_b \|\leq \varepsilon.
\end{eqnarray}
\end{defn}
\begin{prop}\label{Prop1} Let $(H)$ be defined as $(\ref{H})$. Then following properties are equivalent:

$(1)$ The system $(\ref{problem})$ is regionally exactly controllable on $\omega$ at time $b$;\\
$(2)$ $im p_\omega H=L^2(\omega)$;\\
$(3)$ $ ker p_\omega + im H= Z;$\\
$(4)$ For $z\in L^2(\omega),$ there exists a $\gamma>0$ such that
\begin{eqnarray}
\|z\|_{L^2(\omega)}\leq \gamma \|H^*p_\omega^* z\|_{L^2(0,b;\mathbf{R}^p)}.
\end{eqnarray}
\end{prop}
\textbf{Proof.}
Obviously, $(1)\Leftrightarrow (2).$

$(2)\Rightarrow (3):$ For any $z\in L^2(\omega),$ let $\hat{z}$ be the extension of $z$ to $L^2(\Omega).$
Since  $im p_\omega H=L^2(\omega)$, there exists $u\in L^2(0,b;\mathbf{R}^p)$, $z_1\in ker p_\omega$ such that $\hat{z}=z_1+Hu.$

$(3)\Rightarrow (2):$ For any $\tilde{z}\in Z$, from $(3)$, $\tilde{z}=z_1+z_2,$ where $z_1\in ker p_\omega$ and $z_2\in im H$. Then there exists a $u\in L^2(0,b;\mathbf{R}^p)$ such that $Hu=z_2$. Hence, it follows from the definition of $p_\omega$ that $im p_\omega H=L^2(\omega)$.

$(1)\Leftrightarrow (4):$ Here, we note that the equivalence between $(1)$ and $(4)$ can be deduced based on  the following general result in \cite{dualityrelationship}:

Let $E, F, G $ be reflexive Hilbert spaces and $f \in \mathcal{L}(E,G),$ $g \in \mathcal{L}(F,G).$ Then
the following two  properties are equivalent

$(1)$  $imf \subseteq im g$;\\
$(2)$ $\exists ~\gamma > 0$ such that $\|f^*z^*\|_{E^*}\leq \gamma\|g^*z^*\|_{F^*},~~\forall z^*\in G$.\\
By choosing $E=G=L^2(\omega),$ $F=L^2(0,b;\mathbf{R}^p),$ $f=Id_{L^2(\omega)}$ and $g=p_\omega H$, we then obtain the results
and completes the proof.
\begin{prop} There is an equivalence among the following properties:

$\left<1\right>$ The system $(\ref{problem})$ is regionally approximately controllable on $\omega$ at time $b$;\\
$\left<2\right>$ $\overline{im p_\omega H}= L^2(\omega);$\\
$\left<3\right>$ $ ker p_\omega + \overline{im H}= Z;$\\
$\left<4\right>$ The operator $p_\omega HH^*p^*_\omega$ is positive definite.
\end{prop}
\textbf{Proof.}
Similar to the argument in Proposition~\ref{Prop1}, we obtain that $\left<1\right>\Leftrightarrow \left<2\right>\Leftrightarrow \left<3\right>$.
Finally, we show that $\left<2\right>\Leftrightarrow \left<4\right>$.
In fact, it is well known that
\begin{eqnarray*}
\overline{im p_\omega H}= L^2(\omega){\kern 130pt}\\\Leftrightarrow
(p_\omega Hu,z)=0,\forall u\in L^2(0,b;\mathbf{R}^p) \mbox{  implies  } z=0.
\end{eqnarray*}
Let $u=H^*p_\omega^* z$. Then we see that
\begin{eqnarray*}
\overline{im p_\omega H}= L^2(\omega){\kern 126pt} \\\Leftrightarrow
(p_\omega HH^*p_\omega^* z,z)=0 \mbox{  implies  } z=0,~z\in L^2(\omega),
\end{eqnarray*}
i.e.,  the operator $p_\omega HH^*p^*_\omega$ is positive definite and the proof is complete.

\begin{rem}\label{remark7}

$(1)$  The definition~\ref{Def2} can be applied to the case
where $\omega=\Omega$. Note that there exists a system, which is  not controllable on the whole domain
but regionally controllable (see Example 5.1 below).

$(2)$ A system which is exactly (respectively approximately)
controllable on $\omega$ is  exactly (respectively approximately)
controllable on $\omega_1$ for every $\omega_1\subseteq \omega.$
\end{rem}

\section{Regional strategic actuators}\label{sec:3}

In this section, we will explore the characteristic of actuators when the system $(\ref{problem})$ is regionally approximately controllable.

As pointed out in \cite{AEIJAI}, an actuator is a couple $(D,g)$ where
$D\subseteq \Omega $ is the support of the actuator and $g$ is its spatial distribution.
To state our main results, it is supposed that the system under consideration is excited by $p$ actuators $(D_i,g_i)_{1\leq i\leq p}$ and
let $Bu=\sum\limits_{i=1}^p{p_{D_i}g_i(x)u_i(t)}$, where $p\in \mathbf{N}$, $g_i(x)\in L^2(\Omega)$, $u=(u_1,u_2,\cdots,u_p)$ and $u_i(t)\in L^2(0,b)$.
Then the system $(\ref{problem})$ can be rewritten as follows:
\begin{equation}\label{problem3}
\left\{\begin{array}{l}
_0D^{\alpha}_{t}z(t,x)=Az(t,x)+\sum\limits_{i=1}^p{p_{D_i}g_i(x)u_i(t)}\\{\kern 150pt}\mbox{ in }\Omega \times [0,b],
\\\lim\limits_{t\to 0^+} z(t,x)=z_0(x) \mbox{ in }\Omega.
\end{array}\right.
\end{equation}
Moreover, suppose that $-A$ is a uniformly elliptic operator.
By \cite{Hilbert}, we get that
there exists a sequence $(\lambda_j,\xi_{jk}): k=1,2,\cdots,r_j$, $j=1,2,\cdots $
such that

$(1)$ For each $ j=1,2,\cdots$, $\lambda_j$ is the eigenvalue of the operator $-A$ with multiplicities $r_j$ and
\[0< \lambda_1<\lambda_2< \cdots< \lambda_j<\cdots, ~~\lim\limits_{j\to \infty}\lambda_j=\infty.\]

$(2)$ For each $ j=1,2,\cdots$,  $\xi_{jk}~(k=1,2,\cdots,r_j)$ is the orthonormal eigenfunction corresponding to $\lambda_j$, i.e.,
\[(\xi_{jk_m},\xi_{jk_n})=\left\{\begin{array}{l}1,~~k_m=k_n,
\\0,~~k_m\neq k_n,
\end{array}\right.\]
where $1\leq k_m,k_n\leq  r_j,$ $k_m,k_n\in \mathbf{N}$ and  $(\cdot,\cdot)$ is the inner product of space $L^2(\Omega)$.

Then we see that the strongly continuous semigroup $\{\Phi(t)\}_{t\geq 0}$ on  $Z$ generated by $A$ is
\begin{eqnarray}\label{phi}
\Phi(t)z(x)=\sum\limits_{j=1}^{\infty}{\sum\limits_{k=1}^{r_j}{\exp(-\lambda_jt)(z,\xi_{jk})\xi_{jk}(x)}}, ~~x\in \Omega
\end{eqnarray}
and the sequence $\{\xi_{jk}, k=1,2,\cdots,r_j, j=1,2,\cdots \}$ is
an orthonormal basis in $L^2(\Omega)$, then for any $z(x)\in L^2(\Omega)$, it can be expressed as
\begin{eqnarray*}
z(x)=\sum\limits_{j=1}^{\infty}{\sum\limits_{k=1}^{r_j}(z,\xi_{jk})\xi_{jk}(x)}.
\end{eqnarray*}

\begin{defn}
An actuators (or a suite of actuators) is said to be $\omega-$strategic if the system under consideration is regionally approximately controllable on $\omega$.
\end{defn}

Before showing our main result in this part, from Eq.$(\ref{Salpha})$ and Eq.$(\ref{phi})$, for any $z\in L^2(\Omega),$ we have
\begin{eqnarray*}
&~&K_\alpha(t)z(x)\\&=&{\alpha}\int_0^\infty{\theta\phi_{\alpha}(\theta)\Phi(t^{\alpha}\theta)z(x)}d\theta\\
&=&{\alpha}\int_0^\infty{\theta\phi_{\alpha}(\theta)\sum\limits_{j=1}^{\infty}{\sum\limits_{k=1}^{r_j}
{\exp(-\lambda_jt^{\alpha}\theta)(z,\xi_{jk})\xi_{jk}(x)}}}d\theta\\
&=&\sum\limits_{j=1}^{\infty}{\sum\limits_{k=1}^{r_j}\sum\limits_{n=0}^{\infty}\frac{\alpha(-\lambda_jt^{\alpha})^n}{n!}}
(z,\xi_{jk})\xi_{jk}(x)\int_0^\infty{\theta^{n+1}\phi_{\alpha}
}d\theta\\
&=&\sum\limits_{j=1}^{\infty}{\sum\limits_{k=1}^{r_j}\sum\limits_{n=0}^{\infty}\frac{\alpha(n+1)!(-\lambda_jt^{\alpha})^n}{\Gamma(\alpha n+\alpha+1)n!}}
(z,\xi_{jk})\xi_{jk}(x)\\
&=&\sum\limits_{j=1}^{\infty}\sum\limits_{k=1}^{r_j}\alpha E_{\alpha,\alpha+1}^2(-\lambda_jt^{\alpha})(z,\xi_{jk})\xi_{jk}(x),
\end{eqnarray*}
where $E_{\alpha,\beta}^\mu(z):=\sum\limits_{n=0}^{\infty}\frac{(\mu)_n}{\Gamma(\alpha n+\beta)}\frac{z^n}{n!}$, $z\in \mathbf{C}$,
$\alpha,\beta,\mu \in \mathbf{C}$, $\mathbf{Re}{\kern 2pt}\alpha>0$ is the generalized Mittag-Leffler function in three parameters
and here, $(\mu)_n$ is the Pochhammer symbol
defined by (see \cite{H-function}, Section 2.1.1)
\begin{eqnarray}(\mu)_n=\mu(\mu+1)\cdots(\mu+n-1),~n\in \mathbf{N}.\end{eqnarray}
Moreover, If $\alpha,\beta \in \mathbf{C} $  such that $\mathbf{Re}{\kern 2pt}\alpha>0$, $\mathbf{Re}{\kern 2pt} \beta>1$, then
(see Section 2.3.4  in \cite{Mathai-Haubold}, or Section 5.1.1 in \cite{Gorenflo})
\begin{eqnarray}
\alpha E_{\alpha,\beta}^2=E_{\alpha,\beta-1}-(1+\alpha-\beta)E_{\alpha,\beta}.
\end{eqnarray}
It follows that
\begin{eqnarray}\label{3.5}
K_\alpha(t)z(x)
=\sum\limits_{j=1}^{\infty}\sum\limits_{k=1}^{r_j} E_{\alpha,\alpha}(-\lambda_jt^{\alpha})(z,\xi_{jk})\xi_{jk}(x)
\end{eqnarray}
and
\begin{eqnarray*}
&~&\int_0^t\tau^{\alpha-1}{K_\alpha(\tau)Bu(t-\tau)}d\tau\\
&=&\sum\limits_{j=1}^\infty\sum\limits_{k=1}^{r_j}
\sum\limits_{i=1}^p\int_0^t{g^i_{jk}u_i(t-\tau)^{\alpha-1}E_{\alpha,\alpha}
(-\lambda_j\tau^\alpha)}d\tau\xi_{jk}(x),
\end{eqnarray*}
where $E_{\alpha,\beta}(z):=\sum\limits_{i=0}^\infty
{\frac{z^i}{\Gamma(\alpha i+\beta)}},$ $\mathbf{Re}{\kern 2pt}\alpha>0, ~\beta,z\in \mathbf{C}$
 is known as the generalized Mittag-Leffler function in two parameters and $g^i_{jk}=(p_{D_i}g_i,\xi_{jk})$, $j=1,2,\cdots $, $k=1,2,\cdots,r_j$, $i=1,2,\cdots,p$.
Then we have the following result.
\begin{thm} For $j=1,2,\cdots$, arbitrary given $b>0$,
 define $p\times r_j$ matrices $G_j$ as
\begin{equation} \label{G}
G_j=\left[ {\begin{array}{*{20}{c}}
{g_{j1}^1}&{g_{j2}^1}&{\cdots}&{g_{jr_j}^1}\\
{g_{j1}^2}&{g_{j2}^2}&{\cdots}&{g_{jr_j}^2}\\
{\vdots}&{\vdots}&{\vdots}&{\vdots}\\
{g_{j1}^p}&{g_{j2}^p}&{\cdots}&{g_{jr_j}^p}
\end{array}} \right]_{p\times r_j},\end{equation}
where $g^i_{jk}=(p_{D_i}g_i,\xi_{jk})$, $j=1,2,\cdots $, $k=1,2,\cdots,r_j$, $i=1,2,\cdots,p$.
 Then the suite of actuators $(D_i,g_i)_{1\leq i\leq p}$ is said to be $\omega-$strategic  if and only if
\begin{eqnarray}
p\geq r=\max\{r_j\}\mbox{ and }rank ~G_j=r_j,~ j=1,2,\cdots.
\end{eqnarray}
\end{thm}
\textbf{Proof.}
For any given $b>0$ and all $u\in L^2(0,b; \mathbf{R}^p)$, suppose that $z_*\in L^2(\omega)$ satisfies
\begin{eqnarray}\label{3.2.34}
\begin{array}{l}
\left(p_\omega Hu, z_*\right)\\
=\sum\limits_{j=1}^\infty\sum\limits_{k=1}^{r_j}
\sum\limits_{i=1}^p\int_0^b{\frac{E_{\alpha,\alpha}
(-\lambda_j\tau^\alpha)}{\tau^{1-\alpha}}u_i(b-\tau)}d\tau g^i_{jk}z_{jk}
=0,\end{array}\end{eqnarray}
where $z_{jk}=(\xi_{jk},z_*)_{L^2(\omega)}$, $j=1,2,\cdots$, $k=1,2,\cdots,r_j$.
Moreover,  since $u=(u_1,u_2,\cdots, u_p)$ in $(\ref{3.2.34})$ is arbitrary, Lemma~\ref{Lem2} leads us to
\begin{eqnarray}\label{attainaletransform2}
\sum\limits_{j=1}^\infty\sum\limits_{k=1}^{r_j}t^{\alpha-1}E_{\alpha,\alpha}
(-\lambda_jt^\alpha) g^i_{jk}z_{jk}=0.
\end{eqnarray}
Then we see that  the suite of actuators $(D_i,g_i)_{1\leq i\leq p}$ is $\omega-$strategic if and only if
for any $z_*\in L^2(\omega)$, one has
 \begin{eqnarray*}\label{3.10}
\sum\limits_{j=1}^{\infty}{
t^{\alpha-1}E_{\alpha,\alpha}
(-\lambda_jt^\alpha)}G_jz_j=\mathbf{0},~i=1,2,\cdots,p,t\in [0.b]
\end{eqnarray*}
$\Rightarrow z_*=0,$
where $\mathbf{0}=(0,0,\cdots,0)\in \mathbf{R}^p$, $z_j=(z_{j1},z_{j2},\cdots,z_{jr_j})^T$
is a vector in $\mathbf{R}^{r_j}$ and $j=1,2,\cdots$.

Finally, since $t^{\alpha-1}E_{\alpha,\alpha}
(-\lambda_jt^\alpha)>0$ for all $t\geq 0$, $j=1,2,\cdots,$ we then show our proof by using the Reductio and absurdum.

$a)$ If the actuators $(D_i,g_i)_{1\leq i\leq p}$
are not $\omega-$strategic, i.e., the system  $(\ref{problem})$ is not regionally approximately controllable on $\omega$. There exists a
$z_{j^*k}\neq 0$ satisfying \begin{eqnarray}
G_{j^*} z_{j^*} =\mathbf{0}.
\end{eqnarray}
 Then  if
$p\geq r=\max\{r_j\}$, we see that
\begin{eqnarray}
rank ~G_{j^*}<r_{j^* }.
\end{eqnarray}
$b)$ On the contrary, if $p\geq r=\max\{r_j\}$ and $rank ~G_j<r_j \mbox{ for some } j=1,2,\cdots,$ there exists a nonzero
element $\tilde{z} \in L^2( \omega)$ with
$\tilde{z}_j=\left(\tilde{z}_{j1},\tilde{z}_{j2},\cdots, \tilde{z}_{jr_j}\right)^T\in \mathbf{R}^{r_j}$
 such that
\begin{eqnarray}
G_{j} \tilde{z}_{j}=\mathbf{0}.
\end{eqnarray}
Then there exists a nonzero element $\tilde{z} \in  L^2( \omega)$ satisfying
 \begin{eqnarray}
\sum\limits_{j=1}^{\infty}{t^{\alpha-1}E_{\alpha,\alpha}
(-\lambda_jt^\alpha)}G_j\tilde{z}_{j}=\mathbf{0},~t\geq 0.
\end{eqnarray}
This implies that $\overline{im p_\omega H}\neq L^2(\omega)$ and
the suite of actuators $(D_i,g_i)_{1\leq i\leq p}$ is not $\omega-$strategic.
 The proof is complete.

\begin{rem}

$1)$ The system $(\ref{problem})$ with $\alpha=1$,
\[A=-\left(\frac{\partial^2}{\partial x_1^2}+\frac{\partial^2}{\partial x_2^2}+\cdots+\frac{\partial^2}{\partial x_n^2}\right)+q\] and
$q(x)$ being H$\ddot{o}$lder continuous on the compact domain of $\mathbf{R}^n$ is discussed in \cite{179},
 which can be considered as a particular case of our results.

$2)$ If the multiplicity of the eigenvalues $\lambda_j$ of the operator $-A$ is infinite  for some
 $j=1,2,\cdots$ and if the system $(\ref{problem})$ is regionally approximately controllable,
 then  the number of the control functions should not be finite.
\end{rem}

\section{An approach for regional target control}\label{sec:4}

The purpose of this section is to present an approach on how to achieve the regional approximate
 controllability on $\omega$ with the minimum control energy to steer the system  $(\ref{problem})$ from
the initial vector $z_0$ to a target function $z_b$ in the region $\omega$. The method used here
is the Hilbert uniqueness methods (HUMs)\cite{Lions2}.

Let $U_b$ be the closed convex set defined by
\begin{eqnarray}U_b=\{u\in  L^2\left(0,b; \mathbf{R}^p\right): p_\omega z(b,u)=z_b\}.\end{eqnarray}
Consider the following minimization problem
\begin{eqnarray}\label{minimum}
\inf\limits_u J(u)=\inf\limits_u \left\{{\int_0^b{\|u(t)\|^2_{\mathbf{R}^p}}dt}: u\in U_b\right\}. \end{eqnarray}
Next, we show a direct approach to the solution of the regional controllability problem with minimum control energy by
utilizing the HUMs.

Let $G$ and $E$ be the sets given by
\begin{eqnarray}
G=\{g\in L^2(\Omega): g=\mathbf{0}\mbox{ in } \Omega\backslash\omega \}
\end{eqnarray}
and
\begin{eqnarray}
 E=\{e\in L^2(\Omega):e =\mathbf{0}\mbox{ in } \omega \}.
\end{eqnarray}
Then for $(g,e)\in G\times E$, we have
\begin{eqnarray}\label{4.3}
(g,e)=\int_\Omega{ge}dx=\int_\omega{ge}dx+\int_{\Omega\backslash \omega}{ge}dx=0.
\end{eqnarray}
Moreover, for any $g\in G$, consider the system
\begin{equation}\label{Gfunction}
\left\{\begin{array}{l}
Q{}_tD^{\alpha}_b\varphi(t)=A^*Q\varphi(t), ~~ t\in [0,b],
\\ \lim\limits_{t\to 0^+} Q{}_tD^{\alpha-1}_b\varphi(t)=p_\omega^*g
\end{array}\right.
\end{equation}
and the semi-norm on $G$
\begin{equation}\label{Gnorm}
g\in G \to \|g\|_G^2=\int_0^b{\| B^*\varphi (t)\|^2}dt,
\end{equation}
where the reflective operator $Q$ is defined in $(\ref{Q})$.
\begin{lem}\label{Lem4}
$(\ref{Gnorm})$ defines a norm on $G$ if the system $(\ref{problem})$ is  regionally approximately controllable on $\omega$.
\end{lem}
\textbf{Proof.}
For any $g\in G$, by Lemma~\ref{Lem1}, we see that system $(\ref{Gfunction})$ can be rewritten as
\begin{equation}
\left\{\begin{array}{l}
{}_0D^{\alpha}_tQ\varphi(t)=A^*Q\varphi(t), ~~ t\in [0,b],
\\ \lim\limits_{t\to 0^+} {}_0D^{\alpha-1}_tQ\varphi(t)=p_\omega^*g
\end{array}\right.
\end{equation}
and its unique mild solution is
\begin{eqnarray}
\varphi(t)=(b-t)^{\alpha-1}K_\alpha^*(b-t)p_\omega^*g.
\end{eqnarray}
Moreover, if the system $(\ref{problem})$ is regionally approximately controllable on $\omega$,  we have
\begin{eqnarray}\ker  H^* p^*_\omega=\{0\},\end{eqnarray}
i.e.,
\begin{eqnarray}B^*(b-s)^{\alpha-1}K_\alpha^*(b-t)p^*_\omega g=0 \Rightarrow  g=0.\end{eqnarray}
 Hence, for any $g\in G$, it follows from
  \begin{eqnarray*}
&~&\|g\|_G^2=\int_0^b{\|B^*K_\alpha^*(b-s)p^*_\omega g\|^2}ds=0\\
&\Leftrightarrow &\\
&~&B^*(b-s)^{\alpha-1}K_\alpha^*(b-s)p^*_\omega g=\mathbf{0}\end{eqnarray*}
 that    $\|\cdot\|_G$ is a norm of space $G$ and the proof is complete.

In addition, consider the following system
\begin{equation}
\left\{\begin{array}{l}
_0D^{\alpha}_{t}\psi(t)=A\psi(t)+BB^*\varphi(t), ~~ t\in [0,b],
\\ \lim\limits_{t\to 0^+}{}_0D^{\alpha-1}_{t}\psi(t)=0,
\end{array}\right.
\end{equation}
which is controlled by the solution of the system $(\ref{Gfunction})$. Let  $\Lambda_: ~G\to E^\perp$  be
\begin{eqnarray}\label{wedge}
\Lambda_ g=p_\omega \psi(b).
\end{eqnarray}
Suppose that $\widetilde{\psi}(t)$ satisfies
\begin{eqnarray}\label{2.11}
 \left\{\begin{array}{l}
_0D^{\alpha}_{t}\widetilde{\psi}(t)=A\widetilde{\psi}(t),
\\ \lim\limits_{t\to 0^+}{}\widetilde{\psi}(t)=z_0.
\end{array}\right.
\end{eqnarray}
For all $z_b\in L^2(\omega),$  we see that $z_b=p_\omega\left[\psi(b)+ \widetilde{\psi}(b)\right]$ and
the regional controllability  problem is equivalent to solving the equation
\begin{eqnarray}\label{4.13}
\Lambda_ g:= z_b- p_\omega\widetilde{\psi}(b),
\end{eqnarray}
Then we can obtain the following theorem.

\begin{thm}\label{Th2}
If the system $(\ref{problem})$ is  regionally approximately controllable on $\omega$, then for any
$z_b\in L^2(\omega),$ $(\ref{4.13})$ has a unique solution $g\in G$ and the control
\begin{eqnarray}u^*(t)=B^*\varphi(t)\end{eqnarray}
steers the system $(\ref{problem})$ to $z_b$
 at time $b$ in $\omega$. Moreover, $u^*$ solves the minimum problem $(\ref{minimum})$.
\end{thm}
\textbf{Proof.}
By Lemma $\ref{Lem4}$, we see that if the system $(\ref{problem})$ is
regionally approximately controllable on $\omega$, then $\|\cdot\|_G$ is a norm of space $G$.
 Let the completion of $G$ with respect to the norm $\|\cdot\|_G$ again by $G$.
 Then  we will show that $(\ref{4.13})$ has a unique solution in $G$.

For any $g\in G$,
it follows from the definition of operator $\Lambda$ in $(\ref{wedge})$ that
\begin{eqnarray*}\left<g,\Lambda_ g\right>&=&\left<g, p_\omega \psi(b)\right>\\&=&
\left<g,p_\omega\int_0^b(b-s)^{\alpha-1}{K_\alpha(b-s)Bu^*(s)}ds\right>\\&=&
\int_0^b{\left<g,p_\omega (b-s)^{\alpha-1}K_\alpha(b-s)Bu^*(s)\right>}ds \\&=&
\int_0^b{\|B^*\varphi (t)\|^2}ds =\|g\|^2_G.\end{eqnarray*}
Hence, $\Lambda:~ G\to E^\bot$ is one to one. It follows from Theorem 2.1 in \cite{Lions2}
that $(\ref{4.13})$ admits a unique solution in $G$.

Further, let $u=u^*$ in problem $(\ref{problem})$, one has $p_\omega z(b,u^*)=z_b.$
Then for any $u_1\in L^2(0,b,\mathbf{R}^p )$ with $p_\omega z(b,u_1)=z_b$,
we obtain that $p_\omega\left[z(b,u^*)-z(b,u_1)\right]=0.$  Moreover, for any $g\in G,$ we have
$\left<g,p_\omega\left[z(b,u^*)-z(b,u_1)\right]\right>=0$ and
\begin{eqnarray*}
 \begin{array}{l}
0=\left<p_\omega^* g,  \int_0^b(b-s)^{\alpha-1}{K_\alpha(b-s)B(u^*(s)-u_1(s))}ds\right>\\{\kern 5pt}
=\int_0^b{\left<B^*(b-s)^{\alpha-1}K_\alpha^*(b-s)p_\omega^*g,u^*(s)-u_1(s)\right>}ds\\{\kern 5pt}
=\int_0^b{\left<B^*\varphi(t),u^*(s)-u_1(s)\right>}ds.
\end{array}\end{eqnarray*}
By the Theorem 1.3 in \cite{Lions2}, it then follows from
\begin{eqnarray}
\begin{array}{l}
J'(u^*)\cdot(u^*-u_1)=2\int_0^b{\left<u^*(s),u^*(s)-u_1(s)\right>}ds\\{\kern 75pt}
=2\int_0^b{
\left<B^*\varphi(t),u^*(s)-u_1(s)\right>}ds\\{\kern 75pt}
=0,
\end{array}
\end{eqnarray}
that $u^*$ solves the minimum energy problem $(\ref{minimum})$ and the proof is complete.

\section{Examples}\label{sec:5}

This section aims to present two examples to show the effectiveness of our obtained results.

\textbf{Example 5.1.}

 Let us consider the following one dimensional time fractional order differential equations of order
$\alpha\in (0,1)$ with a zone actuator
 to show $(1)$ of Remark $\ref{remark7}$.
\begin{eqnarray}\label{examplel}
\left\{\begin{array}{l}
_0D^{\alpha}_{t}z(x,t)=\frac{\partial^2}{\partial x^2}z(x,t)+p_{[a_1,a_2]}u(t)\mbox{in } [0,1]\times [0,b],
\\ \lim\limits_{t\to 0^+}z(x,t)=z_0(x)\mbox{ in }[0,1],
\\ z(0, t)= z(1, t)=0 \mbox{ in}  [0,b],
\end{array}\right.
\end{eqnarray}
where $Bu=p_{[a_1,a_2]}u$ and
 $0\leq a_1\leq a_2\leq 1$.
 Moreover,  we see that $-A=-\frac{\partial^2}{\partial x^2}$ with
$\lambda_i=i^2\pi^2 ,~\xi_i(x)=\sqrt{2}\sin (i\pi x)$,
$
\Phi(t)z=\sum\limits_{i=1}^{\infty}{\exp(-\lambda_it)(z,\xi_i)_{L^2(0,1)}\xi_i}.
$
and
\begin{eqnarray*}
K_\alpha(t)z(x)=\sum\limits_{i=1}^{\infty}E_{\alpha,\alpha}(-\lambda_it^{\alpha})(z,\xi_i)_{L^2(0,1)}\xi_i(x).
\end{eqnarray*}
Since $A=-\frac{\partial^2}{\partial x^2}$ is a self-adjoint operator, we have
\begin{eqnarray*}
\begin{array}{l}
(H^*z)(t)=\left[B^*(b-t)^{\alpha-1}K_{\alpha}^*(b-t)z\right](t)\\
=B^*(b-t)^{\alpha-1}\sum\limits_{i=1}^{\infty}E_{\alpha,\alpha}(-\lambda_i(b-t)^{\alpha})(z,\xi_i)\xi_i(x)\\
=(b-t)^{\alpha-1}\sum\limits_{i=1}^{\infty}E_{\alpha,\alpha}(-\lambda_i(b-t)^{\alpha})(z,\xi_i)\int_{a_1}^{a_2}{\xi_i(x)}dx.
\end{array}\end{eqnarray*}
By $\int_{a_1}^{a_2}{\xi_i(x)}dx=\frac{\sqrt{2}}{i\pi}
\sin{\frac{i\pi (a_1+a_2)}{2}}\sin{\frac{i\pi (a_2-a_1)}{2}}$, we get that  $Ker (H^*)\neq
\{0\}$ $(\overline{Im (H)}\neq L^2(\omega))$ when $a_2-a_1\in \mathbf{Q}$.
Then the system $(\ref{examplel})$ is not controllable
on $[0,1]$.

Next, we show that there exists a sub-region $\omega \subseteq \Omega$ such that the system $(\ref{examplel})$ is possible regional controllability in $\omega$ at time $b$.

Without loss of generality, let
 $a_1=0,~a_2=1/2$, $z_*=\xi_k, (k=4j,j=1,2,3,\cdots)$. Based on the argument above, $z_*$ is not reachable on $\Omega=[0,1].$
However, since  \[E_{\alpha,\alpha}(t)>0~(t\geq 0)\mbox{ and } \int_0^{1/2}{\xi_i(x)}dx=\frac{\sqrt{2}}{i\pi}\left(1-\cos(i\pi/2)\right),\]
$~~i=1,2,\cdots,$ let $\omega=[1/4,3/4]$, we see that
\begin{eqnarray*}
\begin{array}{l}
(H^*p_\omega^*p_\omega z_*)(t)\\
=\sum\limits_{i=1}^\infty\frac{E_{\alpha,\alpha}(-\lambda_i(b-t)^{\alpha})}{(b-t)^{1-\alpha}}
(\xi_i,\xi_k)_{L^2(\frac{1}{4},\frac{3}{4})}\int_0^{1/2}{\xi_i(x)}dx
\\
=\sum\limits_{i\neq 4j}\frac{\sqrt{2}E_{\alpha,\alpha}(-\lambda_i(b-t)^{\alpha})}{i\pi(b-t)^{1-\alpha}}
\int_{1/4}^{3/4}{\xi_i(x)\xi_{4j}(x)}dx \left[1-\cos(i\pi/2)\right]
\\\neq 0.
\end{array}\end{eqnarray*}
Then $z_*$ is possible regional controllability in $\omega=[1/4,3/4]$ at time $b$.

\textbf{Example 5.2.}

Consider the following time fractional differential equations with a pointwise actuator
\begin{eqnarray}\label{example2}
\left\{\begin{array}{l}
_0D^{\alpha}_{t}z(x,t)=\frac{\partial^2}{\partial x^2}z(x,t)+u(t)\delta(x-\sigma)\mbox{ in }[0,1]\times [0,b],
\\  \lim\limits_{t\to 0^+}z(x,0)=z_0\mbox{ in }[0,1],
\\z(0, t)= z(1, t)=0\mbox{ in} [0,b],
\end{array}\right.
\end{eqnarray}
which is excited by  a pointwise control located at $\sigma\in [0,1].$ Here
 $A=\frac{\partial^2}{\partial x^2}$ generates a strongly continuous semigroup.
In addition, for any $g\in G$, by Lemma $\ref{Lem4}$, we see that
\begin{eqnarray*}g\to \|g\|^2_G
=\int_0^b{\left\|\varphi(s)\right\|^2}ds\end{eqnarray*}
defines a norm on $G$,
where $\varphi$ is the unique mild solution of the following problem
\begin{equation}
\left\{\begin{array}{l}
Q{}_tD^{\alpha}_b\varphi(t)=A^*Q\varphi(t), ~~ t\in [0,b],
\\ \lim\limits_{t\to 0^+} Q{}_tD^{\alpha-1}_b\varphi(t)=p_\omega^*g.
\end{array}\right.
\end{equation}
Now if we consider the following system
\begin{equation}
\left\{\begin{array}{l}
_0D^{\alpha}_{t}\psi(t)=A\psi(t)+\delta(x-\sigma)\varphi(\sigma,t),~~ t\in [0,b],
\\  \lim\limits_{t\to 0^+}{}_0D^{\alpha-1}_{t}\psi(t)=0.
\end{array}\right.
\end{equation}
Let $\omega\subseteq [0,1]$ be a subinterval and let  $\Lambda_: G\to H^\perp$  be
\begin{eqnarray}
\Lambda_ g=p_\omega \psi(b).
\end{eqnarray}
Then the regional controllability  of the example  $(\ref{example2})$ is equivalent to solving the equation
\begin{eqnarray}\label{514}
\Lambda_ g:= z_b- p_\omega\widetilde{\psi}(b),~~\forall z_b\in L^2(\omega)
\end{eqnarray}
where $\widetilde{\psi}(t)$ is the solution of the following system
\begin{eqnarray}
 \left\{\begin{array}{l}
_0D^{\alpha}_{t}\widetilde{\psi}(t)=A\widetilde{\psi}(t),~~ t\in [0,b],
\\  \lim\limits_{t\to 0^+}\widetilde{\psi}(t)=z_0.
\end{array}\right.
\end{eqnarray}
Thus, by Theorem $\ref{Th2}$, we can conclude that
if the example $(\ref{example2})$ is  regionally approximately controllable on some subregion of $[0,1]$, for any
$z_b\in L^2(\omega),$ $(\ref{514})$ admits a unique solution $g\in G$. Moreover,  the control
\begin{eqnarray*}
u^*(t)&=&\varphi(\sigma,t)\\
&=&
\sum\limits_{i=1}^{\infty}(b-t)^{\alpha-1}E_{\alpha,\alpha}(-\lambda_i(b-t)^{\alpha})(p_\omega^*g,\xi_i)
\xi_i(\sigma)\end{eqnarray*}
steers $(\ref{example2})$ to $z_b$ at time $b$ and  $u^*$ solves the minimum control energy problem $(\ref{minimum})$.

\section{CONCLUSIONS}
The purpose of this paper is to investigate the regional controllability of  the Riemann-Liouville
time fractional diffusion equations of order $\alpha \in (0,1)$. The
characterizations of strategic actuators when the control inputs appear in the differential
equations as distributed inputs and an approach on the regional controllability with minimum
energy of the problems $(\ref{problem})$ are solved.
Since $E_1(t)=e^{t}$, $t\geq 0$, together with $(\ref{3.5}),$ we get that our results can be regarded as the
extension of the results in \cite{214} and \cite{179}.

Moreover, the results presented here can also be extended
to  complex fractional order distributed parameter dynamic systems. For instance, the problem of
constrained regional control  of fractional order diffusion systems with more complicated
regional sensing and actuation configurations are of great interest. For more
information on the potential topics related to fractional distributed parameter systems, we
refer the readers to \cite{ge2015JAS} and the references therein.

\begin{ack}
This work was supported  by the Natural Science Foundation of Shanghai (No.15ZR1400800).
\end{ack}



\begin{thebibliography}{29}
\expandafter\ifx\csname natexlab\endcsname\relax\def\natexlab#1{#1}\fi
\providecommand{\url}[1]{\texttt{#1}}
\providecommand{\href}[2]{#2}
\providecommand{\path}[1]{#1}
\providecommand{\DOIprefix}{doi:}
\providecommand{\ArXivprefix}{arXiv:}
\providecommand{\URLprefix}{URL: }
\providecommand{\Pubmedprefix}{pmid:}
\providecommand{\doi}[1]{\href{http://dx.doi.org/#1}{\path{#1}}}
\providecommand{\Pubmed}[1]{\href{pmid:#1}{\path{#1}}}
\providecommand{\bibinfo}[2]{#2}
\ifx\xfnm\relax \def\xfnm[#1]{\unskip,\space#1}\fi

\bibitem[{Chen \& Feng, 2016}]{chen2016adaptive}
\bibinfo{author}{Chen, X.}, \&  \bibinfo{author}{Feng Y.} (\bibinfo{year}{2016}).
\newblock \bibinfo{title}{Adaptive control for continuous-time systems with
  actuator and sensor hysteresis}.
\newblock {\it \bibinfo{journal}{Automatica}}
 {\it \bibinfo{volume}{64}}, \bibinfo{pages}{196--207}.

\bibitem[{Courant \& Hilbert, 1966}]{Hilbert}
\bibinfo{author}{Courant, R.}, \& \bibinfo{author}{Hilbert, D.}
  (\bibinfo{year}{1966}).
\newblock {\it \bibinfo{title}{Methods of mathematical physics}\/}
  volume~\bibinfo{volume}{1}.
\newblock \bibinfo{publisher}{CUP Archive}.

\bibitem[{Dacorogna, 2007}]{Bernard}
\bibinfo{author}{Dacorogna, B.} (\bibinfo{year}{2007}).
\newblock {\it \bibinfo{title}{Direct methods in the calculus of variations,
  {Second} edition}\/} volume~\bibinfo{volume}{78}.
\newblock \bibinfo{publisher}{Springer Science \& Business Media}.

\bibitem[{El~Jai, 1991}]{AEIJAI2}
\bibinfo{author}{El~Jai, A.} (\bibinfo{year}{1991}).
\newblock \bibinfo{title}{Distributed systems analysis via sensors and
  actuators}.
\newblock {\it \bibinfo{journal}{Sensors and Actuators A: Physical}\/},  {\it
  \bibinfo{volume}{29}\/}, \bibinfo{pages}{1--11}.

\bibitem[{El~Jai \& Pritchard, 1988}]{AEIJAI}
\bibinfo{author}{El~Jai, A.}, \& \bibinfo{author}{Pritchard, A.~J.}
  (\bibinfo{year}{1988}).
\newblock {\it \bibinfo{title}{Sensors and controls in the analysis of
  distributed systems}\/}.
\newblock \bibinfo{publisher}{Halsted Press}.

\bibitem[{El~Jai et~al., 1995}]{214}
\bibinfo{author}{El~Jai, A.}, \bibinfo{author}{Simon, M.},
  \bibinfo{author}{Zerrik, E.}, \& \bibinfo{author}{Pritchard, A.}
  (\bibinfo{year}{1995}).
\newblock \bibinfo{title}{Regional controllability of distributed parameter
  systems}.
\newblock {\it \bibinfo{journal}{International Journal of Control}\/},  {\it
  \bibinfo{volume}{62}\/}, \bibinfo{pages}{1351--1365}.

\bibitem[{Erd$\acute{e}$lyi et~al., 1953}]{H-function}
\bibinfo{author}{Erd$\acute{e}$lyi, A.}, \bibinfo{author}{Magnus, W.},
  \bibinfo{author}{Oberhettinger, F.}, \& \bibinfo{author}{Tricomi, F.~G.}  (\bibinfo{year}{1953}).
\newblock \bibinfo{title}{Higher transcendental functions, { Vol.1}}.
\newblock {\it \bibinfo{journal}{McGraw-Hill Book Company, 1953.}\/}

\bibitem[{Fujishiro \& Yamamoto, 2014}]{japankf1}
\bibinfo{author}{Fujishiro, K.}, \& \bibinfo{author}{Yamamoto, M.}
  (\bibinfo{year}{2014}).
\newblock \bibinfo{title}{Approximate controllability for fractional diffusion
  equations by interior control}.
\newblock {\it \bibinfo{journal}{Applicable Analysis}\/},  {\it
  \bibinfo{volume}{93}\/}, \bibinfo{pages}{1793--1810}.

\bibitem[{Ge et~al., 2015}]{ge2015JAS}
\bibinfo{author}{Ge, F.}, \bibinfo{author}{Chen, Y.}, \& \bibinfo{author}{Kou,
  C.} (\bibinfo{year}{2015}).
\newblock \bibinfo{title}{Cyber-physical systems as general distributed
  parameter systems: three types of fractional order models and emerging
  research opportunities}.
\newblock {\it \bibinfo{journal}{IEEE/CAA Journal of Automatica Sinica}\/},  {\it
  \bibinfo{volume}{2}\/}, \bibinfo{pages}{353--357}.

\bibitem[{Ge et~al., 2016{\natexlab{a}}}]{Ge2016IJC}
\bibinfo{author}{Ge, F.}, \bibinfo{author}{Chen, Y.}, \& \bibinfo{author}{Kou,
  C.} (\bibinfo{year}{2016}{\natexlab{a}}).
\newblock \bibinfo{title}{Actuator characterizations to achieve approximate
  controllability for a class of fractional sub-diffusion equations}.
\newblock {\it \bibinfo{journal}{International Journal of Control}\/},
\bibinfo{pages}{1--9}.

\bibitem[{Ge et~al., 2016{\natexlab{b}}}]{Ge2016JMAA}
\bibinfo{author}{Ge, F.}, \bibinfo{author}{Chen, Y.}, \& \bibinfo{author}{Kou,
  C.} (\bibinfo{year}{2016}{\natexlab{b}}).
\newblock \bibinfo{title}{Regional gradient controllability of sub-diffusion
  processes}.
\newblock {\it \bibinfo{journal}{Journal of Mathematical Analysis and Applications}\/},  {\it
  \bibinfo{volume}{440}\/}, \bibinfo{pages}{865--884}.

\bibitem[{Gorenflo et~al., 2014}]{Gorenflo}
\bibinfo{author}{Gorenflo, R.}, \bibinfo{author}{Kilbas, A.~A.},
  \bibinfo{author}{Mainardi, F.}, \& \bibinfo{author}{Rogosin, S.~V.}
  (\bibinfo{year}{2014}).
\newblock {\it \bibinfo{title}{{Mittag-Leffler} functions, related topics and
  applications}\/}.
\newblock \bibinfo{publisher}{Springer}.

\bibitem[{Hilfer, 2000}]{hilfer2000applications}
\bibinfo{author}{Hilfer, R.} (\bibinfo{year}{2000}).
\newblock {\it \bibinfo{title}{Applications of fractional calculus in
  physics}\/} volume \bibinfo{volume}{128}.
\newblock \bibinfo{publisher}{World Scientific}.

\bibitem[{Hilfer \& Anton, 1995}]{hilfer1995fractional}
\bibinfo{author}{Hilfer, R.}, \& \bibinfo{author}{Anton, L.}
  (\bibinfo{year}{1995}).
\newblock \bibinfo{title}{Fractional master equations and fractal time random
  walks}.
\newblock {\it \bibinfo{journal}{Physical Review E}\/},  {\it
  \bibinfo{volume}{51}\/}, \bibinfo{pages}{848--851}.

\bibitem[{Kilbas et~al., 2006}]{Kilbas}
\bibinfo{author}{Kilbas, A.~A.}, \bibinfo{author}{Srivastava, H.~M.}, \&
  \bibinfo{author}{Trujillo, J.~J.} (\bibinfo{year}{2006}).
\newblock {\it \bibinfo{title}{Theory and applications of fractional
  differential equations}\/}.
\newblock \bibinfo{publisher}{Elsevier Science Limited}.

\bibitem[{Lin \& Lu, 2013}]{lin2013laplace}
\bibinfo{author}{Lin, S.-D.}, \& \bibinfo{author}{Lu, C.-H.}
  (\bibinfo{year}{2013}).
\newblock \bibinfo{title}{Laplace transform for solving some families of
  fractional differential equations and its applications}.
\newblock {\it \bibinfo{journal}{Advances in Difference Equations}\/},  {\it
  \bibinfo{volume}{2013}\/}, \bibinfo{pages}{1--9}.

\bibitem[{Lions, 1971}]{Lions2}
\bibinfo{author}{Lions, J.~L.} (\bibinfo{year}{1971}).
\newblock {\it \bibinfo{title}{Optimal control of systems governed by partial
  differential equations}\/} volume \bibinfo{volume}{170}.
\newblock \bibinfo{publisher}{Springer Verlag}.

\bibitem[{Mainardi et~al., 2007}]{mainardi2007sub}
\bibinfo{author}{Mainardi, F.}, \bibinfo{author}{Mura, A.},
  \bibinfo{author}{Pagnini, G.}, \& \bibinfo{author}{Gorenflo, R.}
  (\bibinfo{year}{2007}).
\newblock \bibinfo{title}{Sub-diffusion equations of fractional order and their
  fundamental solutions}.
\newblock In {\it \bibinfo{booktitle}{Mathematical methods in engineering}\/}
  (pp. \bibinfo{pages}{23--55}).
\newblock \bibinfo{publisher}{Springer}.

\bibitem[{Mainardi et~al., 2007}]{Mainardi}
\bibinfo{author}{Mainardi, F.}, \bibinfo{author}{Paradisi, P.}, \&
  \bibinfo{author}{Gorenflo, R.} (\bibinfo{year}{2007}).
\newblock \bibinfo{title}{Probability distributions generated by fractional
  diffusion equations}.
\newblock {\it \bibinfo{journal}{2007.~arXiv preprint arXiv:0704.0320}\/}.

\bibitem[{Klimek, 2009}]{Klimek}
\bibinfo{author}{Klimek, M.} (\bibinfo{year}{2009}).
\newblock {\it \bibinfo{title}{On solutions of linear fractional differential
  equations of a variational type}\/}.
\newblock \bibinfo{publisher}{Publishing Office of Czestochowa University of
  Technology}.

\bibitem[{Mandelbrot, 1983}]{Mandelbrot}
\bibinfo{author}{Mandelbrot, B.~B.} (\bibinfo{year}{1983}).
\newblock {\it \bibinfo{title}{The fractal geometry of nature}\/} volume
  \bibinfo{volume}{173}.
\newblock \bibinfo{publisher}{Macmillan}.

\bibitem[{Mathai \& Haubold, 2008}]{Mathai-Haubold}
\bibinfo{author}{Mathai, A.~M.}, \& \bibinfo{author}{Haubold, H.~J.}
  (\bibinfo{year}{2008}).
\newblock {\it \bibinfo{title}{Special functions for applied scientists}\/}
  volume~\bibinfo{volume}{4}.
\newblock \bibinfo{publisher}{Springer}.

\bibitem[{Metzler \& Klafter, 2000}]{ralf}
\bibinfo{author}{Metzler, R.}, \& \bibinfo{author}{Klafter, J.}
  (\bibinfo{year}{2000}).
\newblock \bibinfo{title}{The random walk's guide to anomalous diffusion: a
  fractional dynamics approach}.
\newblock {\it \bibinfo{journal}{Physics Reports}\/},  {\it
  \bibinfo{volume}{339}\/}, \bibinfo{pages}{1--77}.

\bibitem[{Montroll \& Weiss, 1965}]{Montroll1965}
\bibinfo{author}{Montroll, E.~W.}, \& \bibinfo{author}{Weiss, G.~H.}
  (\bibinfo{year}{1965}).
\newblock \bibinfo{title}{Random walks on lattices. {II}}.
\newblock {\it \bibinfo{journal}{Journal of Mathematical Physics}\/},  {\it
  \bibinfo{volume}{6}\/}, \bibinfo{pages}{167--181}.

\bibitem[{Pritchard \& Wirth, 1978}]{dualityrelationship}
\bibinfo{author}{Pritchard, A.}, \& \bibinfo{author}{Wirth, A.}
  (\bibinfo{year}{1978}).
\newblock \bibinfo{title}{Unbounded control and observation systems and their
  duality}.
\newblock {\it \bibinfo{journal}{SIAM Journal on Control and Optimizatio}\/},  {\it
  \bibinfo{volume}{16}\/}, \bibinfo{pages}{535--545}.

\bibitem[{Renardy \& Rogers, 2006}]{uniformlyelliptic}
\bibinfo{author}{Renardy, M.}, \& \bibinfo{author}{Rogers, R.~C.}
  (\bibinfo{year}{2006}).
\newblock {\it \bibinfo{title}{An introduction to partial differential
  equations}\/} volume~\bibinfo{volume}{13}.
\newblock \bibinfo{publisher}{Springer Science \& Business Media}.

\bibitem[{Sakawa, 1974}]{179}
\bibinfo{author}{Sakawa, Y.} (\bibinfo{year}{1974}).
\newblock \bibinfo{title}{Controllability for partial differential equations of
  parabolic type}.
\newblock {\it \bibinfo{journal}{SIAM Journal on Control and Optimizatio}\/},  {\it
  \bibinfo{volume}{12}\/}, \bibinfo{pages}{389--400}.

\bibitem[{Spears \& Spears, 2012}]{2Spears}
\bibinfo{author}{Spears, W.~M.}, \& \bibinfo{author}{Spears, D.~F.}
  (\bibinfo{year}{2012}).
\newblock {\it \bibinfo{title}{Physicomimetics: Physics-based swarm
  intelligence}\/}.
\newblock \bibinfo{publisher}{Springer Science \& Business Media}.

\bibitem[{Uchaikin \& Sibatov, 2012}]{uchaikin2012fractional}
\bibinfo{author}{Uchaikin, V.}, \& \bibinfo{author}{Sibatov, R.}
  (\bibinfo{year}{2012}).
\newblock \bibinfo{title}{Fractional kinetics in solids: Anomalous charge
  transport in semiconductors}.
\newblock {\it \bibinfo{journal}{Dielectrics and Nanosystems (World Science,
  2013)}\/}, .

\bibitem[{Weinberger, 1962}]{uniformlyelliptic2}
\bibinfo{author}{Weinberger, H.} (\bibinfo{year}{1962}).
\newblock \bibinfo{title}{Symmetrization in uniformly elliptic problems}.
\newblock {\it \bibinfo{journal}{Studies in Mathematical Analysis, Stanford University.
  Press}\/},  (pp. \bibinfo{pages}{424--428}).

\end{thebibliography}
%

\end{document}